\documentclass{amsproc}

\newtheorem{theorem}{Theorem}[section]
\newtheorem{lemma}[theorem]{Lemma}
\newtheorem{prop}[theorem]{Proposition}
\newtheorem{corollary}[theorem]{Corollary}
\newtheorem{conjecture}[theorem]{Conjecture}

\theoremstyle{definition}

\theoremstyle{remark}

\numberwithin{equation}{section}

\newcommand{\C}{{\mathbf C}}

\newcommand{\A}{{\mathcal A}}

\newcommand{\tc}{{\rm\bf {TC}}}

\begin{document}

\title{Topological complexity of generic hyperplane complements}
\author{Sergey Yuzvinsky}
\address{University of Oregon}
\email{yuz@uoregon.edu}
\subjclass{Primary 52C35, 55R80; Secondary 14H10, 98C83}

\keywords{Topological complexity, motion planning algorithm, Schwarz genus, 
hyperplane arrangements.}

\begin{abstract}
We prove that the topological complexity of (a motion planning algorithm on) the complement
of generic complex essential  hyperplane arrangement of $n$ hyperplanes in an
$r$-dimensional linear space is min$\{n+1,2r\}$.
\end{abstract}

\maketitle
\section{Introduction}

In this paper we continue the theme started in \cite{FY} - studying the topological
(motion planning) complexity $\tc(M)$ of the complement $M$ of a
complex hyperplane arrangement.
The number $\tc(X)$ was defined for any path-connected topological 
space $X$ by M.Farber in \cite{F1, F2}. 
This number is of fundamental importance for the 
motion planning problem: $\tc(X)$ determines character of instabilities for all  
motion planning algorithms in $X$.

The main result of this paper can be stated as follows:

\begin{theorem}
\label{main}
Let $M$ be the complement of a complex
central essential arrangement of $n$ hyperplanes in the linear
space $V$ of dimension $r>0$.  Then \hfill\break $\tc(M)={\rm min}\{n+1,2r\}$.
\end{theorem}

\section{The motion planning problem}\label{section1}
\label{definition}

In this section we recall the definitions and results from \cite{F1, F2} that we will
use later in this paper.

Let $X$ be a connected topological space $X$ 
that is homotopy equivalent to a CW complex.
 Let $PX$ be the space of all continuous paths $\gamma: [0,1] \to X$,
equipped with the compact-open topology, and
let $\pi: PX \to X \times X$ be the map assigning the end points to a
path: $\pi (\gamma) = (\gamma(0),
\gamma(1))$. The map $\pi$ is a fibration whose fiber is the
based loop space $\Omega X$.
The {topological complexity of} $X$, denoted by $\tc(X)$, is the smallest
number $k$ such that $X \times X$ can be
covered by open sets $U_1, \dots, U_k$, so that for every $i=1,
\dots, k$ there exists a continuous section $s_i:
U_i \to PX, \pi \circ s_i = 1$.

According to \cite{F2}, {\it a motion planner} in 
$X$ is defined by finitely many subsets $F_1, \dots, F_k$ $\subset X\times X$ and continuous maps
$s_i: F_i\to PX$, where $i=1, \dots, k$, such that:
\begin{enumerate}
\item[(a)] the sets $F_1, \dots, F_k$ are pairwise disjoint (i.e.,$F_i\cap F_j=\emptyset$, 
$i\not= j$), and cover $X\times X$;
\item[(b)] $\pi\circ s_i =1_{F_i}$ for any $i=1, \dots, k$;
\item[(c)] each $F_i$ is an ENR.
\end{enumerate}

The subsets $F_i$ are {\it local domains} 
of the motion planner; the maps $s_i$ are {\it local rules}.

In \cite{F2} it is shown that:
{\it the minimal integer $k$, such that a smooth manifold $X$ admits a motion planner
with $k$ local rules, 
equals $\tc(X)$.}  

The other properties of $\tc(X)$ we will need are:

\smallskip
(i) {$\tc(X)$ depends only on the homotopy type of $X$.}

\smallskip
(ii) {$\tc(X)\leq 2{\rm dim}(X)+1.$}

\smallskip
(iii) {$\tc(X\times S^1)\leq \tc(X)+1.$}

\smallskip
Next result provides a lower bound for $\tc(X)$ in terms of the
cohomology ring  $H^\ast(X)$ with coefficients in a
field. The tensor product $H^\ast(X)\otimes H^\ast(X)$ is also a graded ring
with the multiplication
$$
(u_1\otimes v_1)\cdot (u_2\otimes v_2) = (-1)^{|v_1|\cdot |u_2|}\,
u_1u_2\otimes v_1v_2
$$
where $|v_1|$ and $|u_2|$ are the degrees of the cohomology classes
$v_1$ and $u_2$. The cohomology multiplication
$H^\ast(X)\otimes H^\ast(X)\to H^\ast(X)\label{prod}$ is a ring
homomorphism. Let  $Z \subset H^\ast(X)\otimes
H^\ast(X)$ be the kernel of this homomorphism. The ideal $Z$ is
called {\it the ideal of zero-divisors} of
$H^\ast(X)$. The {\it zero-divisors-cup-length} is the length of the
longest nontrivial product in the ideal of
zero-divisors.  

(iv) {The topological complexity $\tc(X)$ is greater
than the zero-divisors-cup-length of
$H^\ast(X)$.}

The topological complexity $\tc(X)$, as well as 
the Lusternik-Schnirelmann category ${\rm cat}(X)$, are particular cases of the notion of {\it Schwarz
genus} (also known as {\it sectional category}) of a fibration; it was introduced and 
thoroughly studied by A.Schwarz in \cite{Sch}.

\section{Hattori theorem}
\label{Hattori}

In this section we recall the necessary definitions from arrangement theory and
the famous result of Hattori. The details can be found in \cite{OT}.

Let $V$ be a complex linear space of a positive dimension $r$.
An arrangement $\A$ in $V$ is a set
$\{H_1,\ldots,H_n\}$ of $n$ hyperplanes for some $n$. The arrangement is {\it essential}
if $\bigcap_{i=1}^nH_i=0$. In particular for an essential arrangement $n\geq r$.
Fix for each $i$ a functional $\alpha_i\in V^*$ such that
${\rm ker}\alpha_i=H_i$. The arrangement is {\it generic} if for any subset 
$I\subset\overline n=\{1,\ldots,n\}$ with $|I|=r$ the respective set of 
functionals is linearly independent.
In particular each generic arrangement is essential. 

For a generic arrangement $\A$ the homotopy type of the space 
$M=M(\A)=V\setminus\bigcup_{i=1}^nH_i$
is easy to describe. First, in order to give a precise reference we need to reduce
$\A$ to an arrangement of affine hyperplanes. For that choose an element of $\A$, say
$H_n$, put $\overline{H_n}=\{v\in V|\alpha_n(v)=1\}$, 
and put $\overline\A=\{H_i\cap\overline{H_n}|i=1,\ldots,n-1\}$.
The arrangement $\overline\A$ consists of affine hyperplanes in the 
affine space $\overline{H_n}$ of 
dimension $r-1$.
Moreover since $\A$ is generic the affine arrangement $\overline\A$ 
is {\it in general position},
i.e., the intersection of any $p$ hyperplanes from it has codimension $p$ for $p\leq r-1$
and is empty for $p>r-1$. In particular $|\overline\A|=n-1$. Since $M=M(\A)$ is 
the total space of a trivial fiber bundle over 
$\overline M=M(\overline\A)$ with the fiber $\C^*$ we have the homotopy
equivalence $M\approx \overline M\times S^1$ (cf. \cite{OT}, Proposition 5.1).

Now we state Hattori's theorem \cite{OT}, Theorem 5.21.
Denote by $T^m$ the (compact) torus of dimension $m$ and for every
$I\subset\overline m$ put
$$T^m_I=\{(z_1,\ldots,z_m)\in T^m|z_j=1, \ {\rm for}\ j\not\in I\}.$$

\begin{theorem}
\label{hattori}
Let $n>r>1$.
For any general position arrangement of $n-1$ affine hyperplanes in $(r-1)$ -
dimensional space
its complement has the homotopy type of $\overline{M_0}$ where $\overline{M_0}$ is 
the skeleton of dimension $r-1$ of the canonical CW-complex of $T^{n-1}$, i.e.,
$$\overline{M_0}=\bigcup_{|I|=r-1}T^{n-1}_I.$$
\end{theorem}

\begin{corollary}
\label{generic}
For any generic arrangement of $n$ linear hyperplanes in $r$ dimensional space
its complement $M$ has the homotopy type of $M_0$ where
$$M_0=S^1\times\bigcup_{|I|=r-1}T^{n-1}_I.$$
\end{corollary}

\begin{proof}
For $n>r$ it follows immediately from Hattori's theorem. For $n=r$ (in particular for
$r=1$) the arrangement consists of all coordiante hyperplanes whence $M\approx (\C^*)^r
\approx T^r=M_0$.
\end{proof}

The property (i) of $\tc(X)$ allows us 
to focus in the rest of the paper
on calculating $\tc(M_0)$. We will always denote by $n$ the number of hyperplanes in
the generic central arrangement $\A$ we will consider 
and by $r$ the dimension of the ambivalent space $V$.

\section{Low bound}
\label{lowbound}

In this section we use the definition of $M_0$ to describe $H^*(M_0;\C)$ and to exhibit
a low bound on $\tc(M_0)$ using the property (iv).

Denote by $E(n)=\oplus_{i=0}^nE(n)_i$ the exterior 
algebra over $\C$ with $n$ generators of degree one.
Also for every $k$, $0\leq k\leq n$, put  $E(n)^k=E(n)/\oplus_{i>k}E(n)_i$ (a truncated
exterior algebra).

From the description of $M_0$ in Corollary \ref{generic} we have
$$H^*(M_0,\C)=E(1)\otimes E(n-1)^{r-1}$$
where the tensor product is taken in the category of graded algebras.
In particular we have the following lemma.

Denote by $e_0$ a generator of $H^*(S^1)=E(1)$ and by $e_1,\ldots,e_{n-1}$ the generators
of $H^*(\overline{M_0})=E(n-1)^{r-1}$. Also for every 
$I=\{i_1<i_2<\cdots<i_k\}\subset\overline{n-1}$
put $e_I=e_{i_1}\cdots e_{i_k}$.
\begin{lemma}
\label{algebras}
The set
$\{e_0e_I|I\subset\overline{n-1},\ |I|=r-1\}$ is a basis of the
linear space $H^r(M_0,\C)$.
\end{lemma}

Now we define the elements in  
the ideal of zero divisors of $H^*(M_0)\otimes H^*(M_0)$ corresponding to the generators.
Namely put $\overline e_i=1\otimes e_i-e_i\otimes 1$ for every $i=0,1,\ldots,n-1$.

\begin{prop}
\label{product}
Let $k={\rm min}\{n-1,2r-2\}$ and $J\subset\overline{n-1}$ with $|J|=k$. 
Then $\pi=\overline e_0\prod_{i\in J}\overline e_i\not=0$.
\end{prop}
\begin{proof}
The linear space
$H^*(M_0)\otimes H^*(M_0)$ is double graded by the subspaces $H^s\otimes H^t$,
$0\leq s,t\leq r$. It suffices to prove that $(r,k+1-r)$-component $\pi_{r,k+1-r}$
of $\pi$ does not vanish.
Clearly this component is
$$\pi_{r,k+1-r}=\sum_{I\subset J,|I|=r-1}\pm e_0e_I\otimes e_{J\setminus I}.$$
Since $|J\setminus I|=k+1-r\leq r-1$ and $H^*(\overline{M_0})=E(n-1)^{r-1}\subset
H^*(M_0)$ all monomials $e_{J\setminus I}$ belong to a basis of $H^{k+1-r}(M_0)$.
The monomials $e_0e_I$ belong to a basis of $H^{r-1}(M_0)$
by Lemma \ref{algebras}. Hence all the summands of $\pi_{r,k+1-r}$ belong
to a basis of $H^*(M_0)\otimes H^*(M_0)$ whence $\pi_{r,k-r}\not=0$. 
This completes the proof.
\end{proof}

Now the property (iv) of $\tc(X)$ immediately implies the following.

\begin{corollary}
\label{low}

$$\tc(M)=\tc(M_0)\geq {\rm min}\{n+1,2r\}.$$
\end{corollary}

\section{Motion planning}

In this section we prove that the upper bound for $\tc(M_0)$ coincides with the low bound
from the previous section. 

First since $M_0\approx \overline{M_0}\times S^1$ we have by property (iii)

$$\tc(M_0)\leq\tc(\overline{M_0})+1.$$
Now suppose $n+1\geq 2r$. Since ${\rm dim}\overline{M_0}=r-1$ we have using property (ii) 
that $\tc(\overline{M_0})\leq 2r-1$ whence
$$\tc(M_0)\leq 2r={\rm min}\{n+1,2r\}.$$

Thus we have to consider only the case $n+1<2r$.
To find the upper bound in this case we constract an explicit motion
planning for $\overline{M_0}$ with $n$ rules.

\begin{theorem}
\label{planning}
For arbitrary $r\leq n$ there exists a motion planning for $\overline{M_0}$ with $n$ rules.
\end{theorem}
\begin{proof}
First for every $J\subset\overline{n-1}$ we define the close subset $F'_J$  
of $T^{n-1}\times T^{n-1}$ via
$$F'_J=\{(u,u')|u_j=u'_j {\ \rm if\ and\ only\ if\ } j\in J\}$$
and put $F_J=F'_J\cap (\overline{M_0}\times \overline{M_0})$.
Then we put $F_i=\bigcup_{|J|=i}F_J$ for every $i=0,1,\ldots, n-1$.
The sets $F_i$ are pairwise disjoint and cover $\overline{M_0}\times\overline{M_0}$ whence we can
take them as the local domains of the motion planning we are constructing.
Since the sets $F_J$ are also pairwise disjoint 
it suffices now to construct local rules on them, i.e.,
(continuous) sections $s_J:F_J\to P\overline{M_0}$.

For that define an auxiliary function $\tau:S^1\to[0,1]$ by treating $S^1$ (in the rest of
the proof) as the set of all complex numbers of norm 1 and putting
$$\tau(z)=\begin{cases}
\frac{1}{2}(1-\frac{|z-1|}{\sqrt{2}})& \text{if $|z-1|\leq\sqrt{2}$},\\
0&\text{otherwise}.
\end{cases}$$
Notice that $\tau(1)=\frac{1}{2}$.
Also for two points $z\not=z'\in S^1$, $z={\rm exp}[\sqrt{-1}\phi]$,
$z'={\rm exp}[\sqrt{-1}\phi']$, where $0\leq \phi,\phi'< 2\pi$,
define the path $\zeta_{z,z'}$ on $S^1$ via
$\zeta_{z,z'}(t)={\rm exp}[\sqrt{-1}(t\phi+(1-t)\phi')]$ (i.e., the moving with a constant
speed from $z$ to $z'$ along the natural orientation of $\C$).

Now for $(u,u')=((u_1,\ldots,u_{n-1}),(u'_1,\ldots,u'_{n-1}))\in T^{n-1}\times T^{n-1}$
we define $s_J(t)=(s_{J,j}(t))_{j\in\overline{n-1}}$ 
via $s_{J,j}(t)=u_j=u'_j$ for every $t\in[0,1]$ if $j\in J$. If $j\not\in J$ we put
$$s_{J,j}(t)=\begin{cases}
    u_j& \text{if $0\leq t<\tau(u_j)$},\\
    \zeta_{u_j,u'_j}(\frac{t-\tau(u_j)}{1-\tau(u_j)-\tau(u'_J)})
    & \text{if $\tau(u_j)\leq t\leq 1-\tau(u'_j)$},\\
    u'_j& \text{if $1-\tau(u'_j)<t\leq 1$}.
\end{cases}$$

It is clear from the definition that $s_{J}$ is continuous and $s_J(0)=u$, $s_J(1)=u'$.
Also since $\tau$ is continuous and $\zeta_{z,z'}$ depends continously
on $(z,z')$ on $S^1\times S^1$ with the
diagonal deleted we see that $s_J$ is conitiously depending
on $(u,u')$ on $F_J$. It is left to check only that $s_J(t)\in \overline{M_0}$ 
for every $t\in [0,1]$. 
In other words we need to check that for every $t$ we have $s_{J,j}(t)=1$ for at
least $n-r$ values of $j$. 

Suppose that $u\in T^{n-1}_I$ and $u'\in T^{n-1}_{I'}$, $|I|=|I'|=r-1$.
Consider the complements $\overline{I}=\overline{n-1}\setminus I$ and $\overline{I'}=
\overline{n-1}\setminus I'$.
Put $I_0=\overline{I}\cap \overline{I'}$ and fix a bijection $\phi:\overline{I}
\setminus I_0\to \overline{I'}\setminus I_0$
putting $j'=\phi(j)$ for every $j\in \overline{I}\setminus I_0$.
Then if $j\in I_0$ we have $j\in J$ whence $s_{J,j}(t)=u_j=u'_j=1$ for every $t$.
Suppose $j\in \overline{I}\setminus I_0$. Then $\tau(u_j)=\tau(1)=\frac{1}{2}$ whence
$s_{J,j}=u_j=1$ for $t\leq \frac{1}{2}$. On the other hand, $\tau(u'_{j'})=\tau(1)=
\frac{1}{2}$ whence $s_{J,j'}(t)=1$ for $t\geq\frac{1}{2}$. Collecting this data we
see that indeed for arbitrary $t$ there are $n-r$ values of $j$ such that $s_{J,j}(t)=1$
which completes the construction of the motion planning whence also the proof.
\end{proof}

\begin{corollary}
$\tc(\overline{M_0})\leq{\rm min}\{n,2r-1\}$ whence
$\tc(M)=\tc(M_0)\leq{\rm min}\{n+1,2r\}$ and Theorem \ref{main} follows.
\end{corollary}

In all cases where the topological complexity has been computed for hyperplane arrangement
complements it coincides with the low bound given by the zero-divisors-cup-length 
(see property iv in section 2). This justifies the following conjecture.

\begin{conjecture}
For every complex central hyperplane arrangement with the complement $M$
the topological complexity $\tc(M)$ is greater by 1 than
the zero-divisors-cup-length of $H^*(M,\C)$.
\end{conjecture}

\bibliographystyle{amsalpha}

\begin{thebibliography}{7}

\bibitem{F1} M. Farber, \textit{Topological complexity of motion planning},
Discrete Comput. Geom. {\bf 29} (2003), 211-221.

\bibitem{F2} M. Farber, \textit{Instabilities of robot motion}, 
Topology Appl. {\bf 140} (2004), 245-266.

\bibitem{FY} M. Farber, S. Yuzvinsky,
\textit{Topological Robotics: Subspace Arrangements and Collision Free Motion
Planning.} Transl. of AMS {\bf 212}(2004), 145-156.

\bibitem{OT} P. Orlik and H. Terao, Arrangements of hyperplanes, Springer-Verlag, 1992.

\bibitem{Sch} A.  Schwarz. \textit{The genus of a fiber space},
A.M.S. Transl. \textbf{55}(1966), 49 - 140.

\end{thebibliography}

\end{document}